\documentclass[reprint,onecolumn,notitlepage,showpacs,preprintnumbers,amsmath,amssymb,aps,pra]{revtex4-1}

\usepackage{graphicx,color}
\usepackage{dcolumn}
\usepackage{bm}
\usepackage{psfrag}
\usepackage{footnpag}
\usepackage{amsmath,amssymb}
\usepackage{nccfoots}
\usepackage{stmaryrd}

\begin{document}

\preprint{arXiv/}

\title{Closed-form solution of a general three-term recurrence relation}

\author{Ivan~Gonoskov}
\email{ivan.gonoskov@gmail.com}
\affiliation{Department of Physics, Ume\aa\ University, SE-90187 Ume\aa, Sweden}\date{\today}

\begin{abstract}
We present a closed-form solution for n-th term of a general three-term recurrence relation with arbitrary given n-dependent coefficients. The derivation and corresponding proof are based on two approaches, which we develop and describe in detail. First, the recursive-sum theory, which gives the exact solution in a compact finite form using a recursive indexing. Second, the discrete dimensional-convolution procedure, which transforms the solution to the non-recursive expression of n, including a finite number of elementary operations and functions. 

\end{abstract}

\maketitle

\section{Introduction}

	A general three-term recurrence relation is usually defined by the following expression:
\begin{equation}\label{i1}
W_{n+1}=A_{n}W_{n}+B_{n}W_{n-1},
\end{equation}
where $n\geq{1}$ ($n\in{\mathbb{N}}$), $W_{n}$ is unknown function of $n$, $A_{n}$ and $B_{n}$ are arbitrary given functions of $n$, and $W_{0}=~C_{0}, W_{1}=~C_{1}$ are initial conditions (we assume not a trivial case when $C_{1}=C_{0}=0$ and $W_{n}\equiv{0}$).

	This well-known relation has a large number of applications and plays an important role in many areas of mathematics and physics. We note here only few of them to underline the importance of the proposed statement. 

	First, the recurrence relation Eq.(\ref{i1}) corresponds to the finite difference equation for the general second order differential equation with unknown function $f(x)$ and arbitrary given $U(x)$, namely: $f''-U(x)\cdot{f}=0$. Therefore, it is widely used for the analytical and numerical analysis (and approximations) in corresponding physical and mathematical applications, see for example \cite{diff1}.

	Second, three-term recurrence relations appear naturally when one uses Frobenius method for solving some linear differential equations and studying some special functions, see \cite{3t}.
	
	Third, the relation Eq.(\ref{i1}) corresponds to the continued fraction with arbitrary given coefficients. Basically, the three-term recurrence relation corresponds to the expressions for the numerators and denominators of this continued fraction (see chapter IV in \cite{cfr1}), which were derived by Euler. 
	
	In this manuscript we obtain a closed-form solution of a canonical three-term recurrence relation, which is equivalent to Eq.(\ref{i1}) in the case of $(A_{n}\neq{0}, \forall{n})$. Our goal is to obtain the solution consisting of a finite number of terms, rather than a variety of methods with infinite series \cite{infCOD}. In that way we develop two approaches for the expressions with a finite number of terms. We start with the recursive-sum theory, which is presented in sec.III. It allows to obtain the exact solution of the three-term recurrence relation in a compact form using a recursive indexing. Next, we develop the discrete dimensional-convolution procedure, which allows to eliminate recursive indexing and represent the solution as a closed-form expression. It means, that the final closed-form expression depends on arbitrary given coefficients and includes a finite number of elementary operations, such as: ``$\,+\,$'', ``$\,-\,$'', ``$\,\times\,$'', ``$\,\div\,$''; and elementary functions, such as: Heaviside step function (or unit step function) and floor function for integer division.

\section{Preliminaries}	

Let us start from the simplification of the Eq.(\ref{i1}). For this we assume that $A_{n}\neq{0}, \forall{n}$ (opposite case is not usually defined as a standard three-term recurrence relation and should be considered separately). Thus, we can use the following substitution:

\begin{equation}\label{i2}
a_{n+1}=W_{n+1}\cdot{\prod_{i=1}^{n}A_{i}}\;.
\end{equation}

According to this substitution and Eq.(\ref{i1}), $a_{n}$ is fulfilled the following three-term recurrence relation: 

\begin{equation}\label{i3}
a_{n+1}=a_{n}+\frac{B_{n}}{A_{n}A_{n-1}}\cdot{}a_{n-1}\;.
\end{equation}

Since the Eq.(\ref{i1}), Eq.(\ref{i3}) are linear and $n, A_{n}, B_{n}$ are arbitrary, we could assume without loss of generality the following initial conditions: $a_{0}=1, a_{1}=1$. Finally, we specify $d_{n}=B_{n}\cdot{}(A_{n}A_{n-1})^{-1}$ as an arbitrary given function of n, and write below the expression, which we further call \textit{canonical three-term recurrence relation}:
\begin{equation}\label{i4}
\begin{aligned}
&a_{n+1}=a_{n}+d_{n}\cdot{}a_{n-1} \;, \\
&a_{0}=1, a_{1}=1.
\end{aligned}
\end{equation}

	We note, that the solution of Eq.(\ref{i4}) consist of a sum of $d_{i}$-products ($i\in{\mathbb{N}}$) with different \textit{powers}. Here we mean that the product \textit{power} $p$ is equal to the total number of $d_{i}$ in the product. The maximum power of $a_{n+1}$, i.e. the largest power among all the products in $a_{n+1}$, is equal to $\llbracket(n+1)/2\rrbracket$, since:

\begin{equation}\label{i5}
\llbracket(n+1)/2\rrbracket=\max\Big(\;\llbracket{}n/2\rrbracket\; , \;1+\llbracket(n-1)/2\rrbracket\;\Big), \forall{n},
\end{equation}
where we use Gauss notation for the integer division: $\llbracket{}x\rrbracket=\left\lfloor x \right\rfloor=\max\big(m\in{\mathbb{Z}} \; | \; m\leq{}x\big)$. Because of this fundamental characteristic it is natural to express $a_{n+1}$, i.e. the solution of Eq.(\ref{i4}), in the following form:

\begin{equation}\label{i6}
a_{n+1} = 1 + \sum_{p=1}^{\llbracket(n+1)/2\rrbracket}S(n,p),
\end{equation}
where $S(n,p)$ includes all terms of the \textit{power} $p$. In the next sections we demonstrate that each $S(n,p)$ is equal to the corresponding recursive-sum, and rigorously proof all the propositions.

\section{recursive-sum theory}

In this section we start from the definition of a general recursive-sum (R-sum). Next, we present and proof some of its properties, which could determine the R-sum algebra. It could be useful for solving different recurrence relations, but, particularly for the solution of the \textit{canonical three-term recurrence relation} it is sufficient to use a particular case of R-sum, namely reduced R-sum which is described in the next section.  
	
	\underline{\textbf{Definition I}}: A general R-sum ($\mathbf{R}$) is defined for the ordered sequence $(d_{i} , i\in{\mathbb{Z}})$ by the following expression:

\begin{equation}\label{i7}
\mathbf{R}=\mathbf{R}(N,k,\Delta,\Delta_{0})=\displaystyle\prod_{m=1}^{k} \left[ \displaystyle \sum_{i_{m}=i_{m-1}+\Delta}^{N+m\cdot\Delta} d_{i_{m}} \right] \;,
\end{equation}
where $(N,k,\Delta,\Delta_{0} \in {\mathbb{Z}})$, $i_{0}\equiv\Delta_{0}$, $N$ is an order of the R-sum, $k\geq{1}$ is a power of the R-sum, $\Delta$ is a recursive shift of the R-sum, $\Delta_{0}$ is an initial shift of the R-sum. For the correct definition it is also assumed that $N\geq\Delta_{0}$ and $d_{i}$ is defined for all indexes used in Eq.(\ref{i7}). The product in Eq.(\ref{i7}) corresponds to the standard left-to-right order of multiplication, so the definition could be written also in the following way:

\begin{equation}\label{i8}
\begin{aligned}
&\mathbf{R}(N,1,\Delta,\Delta_{0})=\displaystyle \sum_{i_{1}=\Delta_{0}+\Delta}^{N+\Delta} d_{i_{1}} \;, &&&&&&&&&&&&&&&&&&&&&&&& \\
&\mathbf{R}(N,2,\Delta,\Delta_{0})=\displaystyle \sum_{i_{1}=\Delta_{0}+\Delta}^{N+\Delta} d_{i_{1}} \cdot{} \displaystyle \sum_{i_{2}=i_{1}+\Delta}^{N+2\Delta} d_{i_{2}} \;, \\
&\mathbf{R}(N,3,\Delta,\Delta_{0})=\displaystyle \sum_{i_{1}=\Delta_{0}+\Delta}^{N+\Delta} d_{i_{1}} \cdot{} \displaystyle \sum_{i_{2}=i_{1}+\Delta}^{N+2\Delta} d_{i_{2}} \cdot{} \displaystyle \sum_{i_{3}=i_{2}+\Delta}^{N+3\Delta} d_{i_{3}} \;, \\
&\;\;\;\;\;\;\;\;\;\,\vdots \quad\quad\quad\quad\quad\quad\quad\quad\quad\quad\quad\quad\;\;\;\;\;\;\;\;\;\;\;\;\;\;\; \cdots \; \\
&\mathbf{R}(N,k,\Delta,\Delta_{0})=\displaystyle \sum_{i_{1}=\Delta_{0}+\Delta}^{N+\Delta} d_{i_{1}} \cdot{} \displaystyle \sum_{i_{2}=i_{1}+\Delta}^{N+2\Delta} d_{i_{2}}  \cdot{}\,\cdots{}\,\cdot{} \displaystyle \sum_{i_{k-1}=i_{k-2}+\Delta}^{N+(k-1)\Delta} d_{i_{k-1}} \cdot{} \displaystyle \sum_{i_{k}=i_{k-1}+\Delta}^{N+k\Delta} d_{i_{k}}  \;. \\
\end{aligned}
\end{equation} 

\newpage
\begin{flushleft}
{Next, we present some of the R-sum elementary properties P.(1-3), which are following directly from the definition.}
\end{flushleft}
\underline{P.(1)}:	
\begin{equation}\label{i9}
\mathbf{R}(N,k,\Delta,N)= \displaystyle \prod_{j=1}^{k} d_{N+j\cdot{}\Delta} \;.
\end{equation}
\underline{P.(2)}:
\begin{equation}\label{i10}
\mathbf{R}(N,k+1,\Delta,\Delta_{0})=\displaystyle \sum_{j=\Delta_{0}+\Delta}^{N+\Delta} d_{j} \cdot{} \mathbf{R}(N+\Delta,k,\Delta,j) \;.
\end{equation}
\underline{P.(3)}:
\begin{equation}\label{i11}
\mathbf{R}(N,k,\Delta,\Delta_{0}+1)= \mathbf{R}(N,k,\Delta,\Delta_{0}) - d_{\Delta_{0}+\Delta}\cdot{}\mathbf{R}(N+\Delta,k-1,\Delta,\Delta_{0}+\Delta) \;.
\end{equation}
\vspace{0.1cm}
\begin{flushleft}{Finally, we consider the \underline{R-sum \textbf{key lemma}}: If $N\geq\Delta_{0}+1$, then:}\end{flushleft}
\begin{equation}\label{i12}
\vspace{0.05cm}
\mathbf{R}(N,k+1,\Delta,\Delta_{0})= \mathbf{R}(N-1,k+1,\Delta,\Delta_{0}) + d_{N+(k+1)\Delta}\cdot{}\mathbf{R}(N,k,\Delta,\Delta_{0}) \;.
\end{equation}
\begin{flushleft}{\underline{\textbf{Proof of the key lemma}}:}\end{flushleft} First, we consider the difference $\mathbf{D}=\left[\mathbf{R}(N,k+1,\Delta,\Delta_{0})-d_{N+(k+1)\Delta}\cdot{}\mathbf{R}(N,k,\Delta,\Delta_{0})\right]$ and try to eliminate index numbers which correspond to zero terms in the original expression. For the correct procedure (when an upper bound in each sum is always not less than the lower one) we will eliminate zero terms and delete the corresponding index numbers step by step from $i_{1}$ to $i_{k+1}$. 
\begin{equation}\label{i13}
\mathbf{D}= \displaystyle \sum_{i_{1}=\Delta_{0}+\Delta}^{N+\Delta} d_{i_{1}} \cdot{} \displaystyle \sum_{i_{2}=i_{1}+\Delta}^{N+2\Delta} d_{i_{2}} \cdot{}\,\cdots{}\,\cdot{} \displaystyle \sum_{i_{k}=i_{k-1}+\Delta}^{N+k\Delta} d_{i_{k}} \cdot{} \left[ \displaystyle \sum_{i_{k+1}=i_{k}+\Delta}^{N+(k+1)\Delta} d_{i_{k+1}}  - d_{N+(k+1)\Delta}\right]\;.
\end{equation}
At the next step we note that for $i_{1}=N+\Delta$ according to the definition we have only one value $i_{2}=N+2\Delta$ and one term $d_{N+2\Delta}$ in the second sum in Eq.(\ref{i13}). According to the definition and P.(1) the same is true for the other indexes: $i_{3}=N+3\Delta$, ..., $i_{k}=N+k\Delta$. The last term in Eq.(\ref{i13}) gives $\left[ \displaystyle \sum_{i_{k+1}=N+k\Delta+\Delta}^{N+(k+1)\Delta} d_{i_{k+1}}  - d_{N+(k+1)\Delta}\right]\equiv{0}$. This means that all the terms in Eq.(\ref{i13}) which correspond to $i_{1}=N+\Delta$ are equal to zero. Since the condition ($N\geq\Delta_{0}+1$), we could then write:
\begin{equation}\label{i14}
\mathbf{D}= \displaystyle \sum_{i_{1}=\Delta_{0}+\Delta}^{(N-1)+\Delta} d_{i_{1}} \cdot{} \displaystyle \sum_{i_{2}=i_{1}+\Delta}^{N+2\Delta} d_{i_{2}} \cdot{}\,\cdots{}\,\cdot{} \displaystyle \sum_{i_{k}=i_{k-1}+\Delta}^{N+k\Delta} d_{i_{k}} \cdot{} \left[ \displaystyle \sum_{i_{k+1}=i_{k}+\Delta}^{N+(k+1)\Delta} d_{i_{k+1}}  - d_{N+(k+1)\Delta}\right]\;.
\end{equation}
The same procedure can be repeated for any successive indexes. Consider the $m$-th sum ($m\leq{}k$) in Eq.(\ref{i13}). For each $i_{m} = N+m\Delta$ we have only one term in the next sum, corresponding to $i_{m+1} = N+(m+1)\Delta$, and so on until the last term: $\left[ \displaystyle \sum_{i_{k+1}=N+k\Delta+\Delta}^{N+(k+1)\Delta} d_{i_{k+1}}  - d_{N+(k+1)\Delta}\right]\equiv{0}$. So we can delete all the corresponding index numbers and write:
\begin{equation}\label{i15}
\mathbf{D}= \displaystyle \sum_{i_{1}=\Delta_{0}+\Delta}^{(N-1)+\Delta} d_{i_{1}} \cdot{} \displaystyle \sum_{i_{2}=i_{1}+\Delta}^{(N-1)+2\Delta} d_{i_{2}} \cdot{}\,\cdots{}\,\cdot{} \displaystyle \sum_{i_{k}=i_{k-1}+\Delta}^{(N-1)+k\Delta} d_{i_{k}} \cdot{} \left[ \displaystyle \sum_{i_{k+1}=i_{k}+\Delta}^{N+(k+1)\Delta} d_{i_{k+1}}  - d_{N+(k+1)\Delta}\right]\;.
\end{equation} 
Finally, we note that $\left[ \displaystyle \sum_{i_{k+1}=i_{k}+\Delta}^{N+(k+1)\Delta} d_{i_{k+1}}  - d_{N+(k+1)\Delta}\right]=\displaystyle \sum_{i_{k+1}=i_{k}+\Delta}^{(N-1)+(k+1)\Delta} d_{i_{k+1}}$. Thus we obtain:
\begin{equation}\label{i16}
\mathbf{D}= \displaystyle \sum_{i_{1}=\Delta_{0}+\Delta}^{(N-1)+\Delta} d_{i_{1}} \cdot{} \displaystyle \sum_{i_{2}=i_{1}+\Delta}^{(N-1)+2\Delta} d_{i_{2}} \cdot{}\,\cdots{}\,\cdot{} \displaystyle \sum_{i_{k}=i_{k-1}+\Delta}^{(N-1)+k\Delta} d_{i_{k}} \cdot{} \displaystyle \sum_{i_{k+1}=i_{k}+\Delta}^{(N-1)+(k+1)\Delta} d_{i_{k+1}} = \mathbf{R}(N-1,k+1,\Delta,\Delta_{0}) \;,
\end{equation} 
which proves the key lemma.	

\newpage

\section{exact solution of a canonical three-term recurrence relation via finite R-sum expansion}

In this section we start with a definition of a reduced R-sum, a particular case of the general R-sum. Then, we construct the exact solution of the \textit{canonical three-term recurrence relation}, by using the R-sum key lemma.  

\underline{\textbf{Definition II}}: A reduced R-sum ($\mathbf{\tilde{R}}$) is a particular case of the general R-sum (see Eqs.(\ref{i7}-\ref{i8})), and defined by the following expression:

	\begin{equation}\label{i17}
 \mathbf{\tilde{R}}=
 \mathbf{\tilde{R}}(N,k)=\mathbf{R}(N,k,2,-1)=\displaystyle \sum_{i_{1}=1}^{N+2} d_{i_{1}} \cdot{} \displaystyle \sum_{i_{2}=i_{1}+2}^{N+4} d_{i_{2}} \cdot{}\;\cdots{}\;\cdot{} \displaystyle \sum_{i_{k-1}=i_{k-2}+2}^{N+2(k-1)} d_{i_{k-1}} \cdot{} \displaystyle \sum_{i_{k}=i_{k-1}+2}^{N+2k} d_{i_{k}} \;.
	\end{equation}
The key lemma Eq.(\ref{i12}) gives the following for the reduced R-sum:
	\begin{equation}\label{i18}
\mathbf{\tilde{R}}(N,k+1) = \mathbf{\tilde{R}}(N-1,k+1) + d_{N+2(k+1)}\cdot{}\mathbf{\tilde{R}}(N,k) \;.
	\end{equation}
Since $N$ and $k$ are independent and arbitrary, let us consider them as functions of new numbers $n$ and $p$: $N=~N(n,p)=n-2p\,;\; k=~k(n,p)=p\,$. In that way we consider a new subsidiary function:
\begin{equation}\label{i19}
S(n,p) = \mathbf{\tilde{R}}(n-2p,p) \;.
\end{equation}
According to the previous notations we have: $N(n,p)=n-2p=n-2k=N(n+2,p+1)=N(n+1,p+1)+1$. Thus we can write the analog of the key lemma for the $S$ function, by using Eq.(\ref{i18}) and Eq.(\ref{i19}):
\begin{equation}\label{i20}
S(n+2,p+1) = S(n+1,p+1) + d_{n+2}\cdot{}S(n,p) \;.
\end{equation}
\vspace{0.1cm}
Now, we will rigorously proof a theorem about the exact solution of a \textit{canonical three-term recurrence relation} Eq.(\ref{i4}).
\begin{flushleft}{\underline{\textbf{Theorem I}}}\end{flushleft} Exact solution of Eq.(\ref{i4}) is following:
\begin{equation}\label{i21}
a_{n+1} = 1 + \displaystyle \sum_{p=1}^{\llbracket(n+1)/2\rrbracket} S(n,p) \;.  
\end{equation}
\begin{flushleft}{\underline{Proof of the \textbf{theorem I.}}}\end{flushleft}
Our proof is based on mathematical induction. Base cases:
\begin{equation}\label{i22}
\begin{aligned}
&n=1:\; &a_{2}=1+d_{1}=1+\displaystyle \sum_{p=1}^{1}S(1,p).\\
&n=2:\; &a_{3}=1+d_{1}+d_{2}=1+\displaystyle \sum_{p=1}^{1}S(2,p),\\
&n=3:\; &a_{4}=1+d_{1}+d_{2}+d_{3}+d_{1}d_{3}=1+\displaystyle \sum_{p=1}^{2}S(3,p),\\
&n=4:\; &a_{5}=1+d_{1}+d_{2}+d_{3}+d_{4}+d_{1}d_{3}+d_{1}d_{4}+d_{2}d_{4}=1+\displaystyle \sum_{p=1}^{2}S(4,p).
\end{aligned}
\end{equation}
Inductive step: we will proof that the statement Eq.(\ref{i21}) for $a_{n}$ and $a_{n-1}$, $\forall{n}$ gives the following:
\begin{equation}\label{i23}
a_{n} + d_{n}\cdot{}a_{n-1} = a_{n+1}.
\end{equation}
To proof it for all $n$, we consider below two cases. First case: $n$ is an arbitrary even number, $n=2m, m\in{\mathbb{N}}$. Second case: $n$ is an arbitrary odd number, $n=2m+1, m\in{\mathbb{N}}$.\\
\underline{First case, $n=2m$.}
\begin{equation}\label{i24}
a_{n} + d_{n}\cdot{}a_{n-1} = a_{2m} + d_{2m}\cdot{}a_{2m-1}=1+d_{2m}+\displaystyle \sum_{l=1}^{m} S(2m-1,l) + d_{2m}\cdot{}\displaystyle \sum_{p=1}^{m-1} S(2m-2,p).
\end{equation}
Note, that $\displaystyle \sum_{l=1}^{m} S(2m-1,l) = S(2m-1,1)+\displaystyle \sum_{l=2}^{m} S(2m-1,l) = S(2m-1,1)+\displaystyle \sum_{p=1}^{m-1} S(2m-1,p+1)$. After the substitution we have:
\begin{equation}\label{i25}
\begin{aligned}
a_{n} + d_{n}\cdot{}a_{n-1} &=1+d_{2m}+ S(2m-1,1)+\displaystyle \sum_{p=1}^{m-1} S(2m-1,p+1) + d_{2m}\cdot{}\displaystyle \sum_{p=1}^{m-1} S(2m-2,p) \\
&=1+d_{2m}+ S(2m-1,1) + \displaystyle \sum_{p=1}^{m-1} \Big[ S(2m-1,p+1) + d_{2m}\cdot{}S(2m-2,p) \Big].
\end{aligned}
\end{equation}
According to Eq.(\ref{i20}) we could write: $S(2m-1,p+1) + d_{2m}\cdot{}S(2m-2,p) = S(2m,p+1)$. Also, according to the definitions Eq.(\ref{i17}) and Eq.(\ref{i19}) we could write: $d_{2m}+ S(2m-1,1) = d_{2m}+ \displaystyle \sum_{j=1}^{2m-1} d_{j} = \displaystyle \sum_{j=1}^{2m} d_{j} = S(2m,1)$. Thus, we could obtain:
\begin{equation}\label{i26}
\begin{aligned}
a_{n} + d_{n}\cdot{}a_{n-1} &=1+S(2m,1)+\displaystyle \sum_{p=1}^{m-1} S(2m,p+1) \\
&= 1+S(2m,1)+\displaystyle \sum_{p=2}^{m} S(2m,p) = 1 + \displaystyle \sum_{p=1}^{m} S(2m,p) = a_{n+1},
\end{aligned}
\end{equation}
which proves the first case of the theorem.\\
\underline{Second case, $n=2m+1$.}
\begin{equation}\label{i27}
a_{n} + d_{n}\cdot{}a_{n-1} = a_{2m+1} + d_{2m+1}\cdot{}a_{2m}=1+d_{2m+1}+\displaystyle \sum_{l=1}^{m} S(2m,l) + d_{2m+1}\cdot{}\displaystyle \sum_{p=1}^{m} S(2m-1,p).
\end{equation}
Similar to the previous case, we write: $\displaystyle \sum_{l=1}^{m} S(2m,l)=S(2m,1)+\displaystyle \sum_{p=1}^{m-1} S(2m,p+1)$. Also we note that $\displaystyle \sum_{p=1}^{m} S(2m-1,p) = \displaystyle \sum_{p=1}^{m-1} S(2m-1,p) + S(2m-1,m)$. Thus, we could write:
\begin{equation}\label{i28}
a_{n} + d_{n}\cdot{}a_{n-1} = 1+d_{2m+1}+S(2m,1)+d_{2m+1}S(2m-1,m)+\displaystyle \sum_{p=1}^{m-1}\Big[ S(2m,p+1) + d_{2m+1}\cdot{}S(2m-1,p) \Big].
\end{equation}
According to Eq.(\ref{i20}) we again could rewrite: $S(2m,p+1) + d_{2m+1}\cdot{}S(2m-1,p) = S(2m+1,p+1)$. In addition we note that $\displaystyle \sum_{p=1}^{m-1}S(2m+1,p+1) = \displaystyle \sum_{p=2}^{m}S(2m+1,p)$. Summarizing, we could obtain:
\begin{equation}\label{i29}
a_{n} + d_{n}\cdot{}a_{n-1} = 1+d_{2m+1}+S(2m,1)+d_{2m+1}S(2m-1,m)+\displaystyle \sum_{p=2}^{m}S(2m+1,p).
\end{equation}
Finally, we note that $d_{2m+1}S(2m-1,m) = d_{2m+1}\mathbf{\tilde{R}}(-1,m) = d_{2m+1}\displaystyle \prod_{j=1}^{m}d_{2j-1}=\mathbf{\tilde{R}}(-1,m+1)=S(2m+1,m+1)$. In addition, we could write $d_{2m+1}+S(2m,1) = d_{2m+1} + \displaystyle \sum_{j=1}^{2m}d_{j} = \displaystyle \sum_{j=1}^{2m+1}d_{j} = S(2m+1,1)$. This gives us the final result:
\begin{equation}\label{i30}
a_{n} + d_{n}\cdot{}a_{n-1} = 1+S(2m+1,1)+S(2m+1,m+1)+\displaystyle \sum_{p=2}^{m}S(2m+1,p)=1+\displaystyle \sum_{p=1}^{m+1}S(2m+1,p) = a_{n+1},
\end{equation}
which proves the second case and the theorem I.

\newpage

\section{discrete dimensional-convolution procedure}

In this section we develop a procedure which could be used for the transformation of the exact solution Eq.(\ref{i21}) to the expression without recursive indexing.

Let us start from a simple case to demonstrate the idea of the discrete dimensional-convolution procedure. We consider a simple expression, which could be associated with a particular case of R-sum.
\begin{equation}\label{i31}
E_{2} = \displaystyle \sum_{i_{1}=1}^{N_{1}}\;\displaystyle \sum_{i_{2}=i_{1}+\delta}^{N_{2}} \, f(i_{1},i_{2}) \,,
\end{equation}
where $N_{1},N_{2}\in{\mathbb{N}}$, $(N_{1}\leq{}N_{2})$ are arbitrary given natural numbers; $\delta\in{\mathbb{N}}$, $(0\leq{}\delta\leq{}[N_{2}-N_{1}])$ is an arbitrary given natural number (shift); and $f(i_{1},i_{2})$ is an arbitrary given function of natural index numbers $i_{1}$ and $i_{2}$. The previous expression $E_{2}$ could be rewritten as follows:
\begin{equation}\label{i32}
E_{2} = \displaystyle \sum_{i_{1}=1}^{N_{1}}\displaystyle \sum_{i_{2}=1}^{N_{2}} F(i_{1},i_{2}) \,,
\end{equation}
where $F(i_{1},i_{2}) = f(i_{1},i_{2})\cdot{}H(i_{2}-i_{1}-\delta)$, and $H(x)$ is Heaviside step function (or unit step function):
\begin{equation}\label{i33}
H(x) = \left\{  \begin{array}{ll}  1\;, & x\geq{0}\\  0\;, & x<0 \end{array}  \right .
\end{equation}
Let us consider now the expression Eq.(\ref{i32}). The indexes $i_{1}$ and $i_{2}$ could be associated with two dimensions by the following way. Let us consider a two dimensional plot (array with elements $F(i_{1},i_{2})$) with horizontal numbering related to $i_{1}$ index and vertical numbering related to $i_{2}$ index, see the figure below.
\begin{equation*}
\begin{bmatrix}
F(1,1) & F(1,2) & \;\cdots\; & F(1,N_{1}) \\
F(2,1) & F(2,2) & \;\cdots\; & F(2,N_{1}) \\
\vdots &        & \;\ddots\; & \vdots \\
F(N_{2},1) & F(N_{2},2) & \;\cdots\; & F(N_{2},N_{1})
\end{bmatrix}
\end{equation*}
For numbering of this array of elements it is possible to use a traversal rule with one global index $q\in{[1,\,\ldots\,,N_{1}\cdot{}N_{2}]}$ instead of $i_{1}\in{[1,\,\ldots\,,N_{1}]}$ and $i_{2}\in{[1,\,\ldots\,,N_{2}]}$. Consider, for example, the following rule: $q = i_{1} + N_{1}\cdot{}(i_{2}-1)$, with the corresponding relations\Footnote{*}{We use the following notations: $\big\{x\big\} = x - \llbracket{}x\rrbracket, \; \llbracket{}x\rrbracket=\left\lfloor x \right\rfloor=\max\big(m\in{\mathbb{Z}} \; | \; m\leq{}x\big)$.}: $i_{1}=i_{1}(q)=1+N_{1}\big\{(q-1)/N_{1}\big\}=q-N_{1}\llbracket(q-1)/N_{1}\rrbracket$, $i_{2}=i_{2}(q)=1+\llbracket(q-1)/N_{1}\rrbracket$. Since, for the presented rule, a unique combination of $(i_{1}, i_{2})$ corresponds to a certain unique number $q$ and vice versa (i.e. we have one-to-one mapping), we could rewrite the Eq.(\ref{i32}) in the following form:
\begin{equation}\label{i34}
E_{2} = \displaystyle \sum_{q=1}^{N_{1}\cdot{}N_{2}}\,F\Big(1+N_{1}\big\{(q-1)/N_{1}\big\} \,,\, 1+\llbracket(q-1)/N_{1}\rrbracket\Big) \,.
\end{equation}
The presented procedure, i.e. the reduction from two dimensions to one dimension, is a particular case of \textit{discrete dimensional-convolution procedure}. Let us now generalize it to a multidimensional case. In that way we consider a multidimensional analog of Eq.(\ref{i32}):
\begin{equation}\label{i35}
E_{k} = \displaystyle \sum_{i_{1}=1}^{N_{1}}\displaystyle \sum_{i_{2}=1}^{N_{2}} \cdots \displaystyle \sum_{i_{k}=1}^{N_{k}} F(i_{1},i_{2}, \ldots{} , i_{k}) \,,
\end{equation}
where ($k\in{\mathbb{N}}$), $k\geq{3}$ - is an arbitrary given natural number, $(N_{1},N_{2},\ldots{},N_{k}\,\in{\mathbb{N}})$ - is an array of arbitrary given natural numbers, and $F(i_{1},i_{2}, \ldots{} , i_{k})$ - is an arbitrary given function of $i_{1}$, $i_{2}$, $\ldots{}$, $i_{k}$.

The previous expression is related to the calculation of a reduced R-sum. According to the definition of the reduced R-sum Eq.(\ref{i17}) and the substitution from Eqs.(\ref{i31}-\ref{i32}) we could write:
\begin{equation}\label{i36}
\mathbf{\tilde{R}}(N,k) = \displaystyle \sum_{i_{1}=1}^{N_{1}}\displaystyle \sum_{i_{2}=1}^{N_{2}} \cdots \displaystyle \sum_{i_{k}=1}^{N_{k}} \Bigg( d_{i_{1}} \cdot{} \displaystyle\prod_{m=2}^{k} \Big[d_{i_{m}}H(i_{m}-i_{m-1}-2)\Big] \Bigg) = \displaystyle \sum_{i_{1}=1}^{N_{1}}\displaystyle \sum_{i_{2}=1}^{N_{2}} \cdots \displaystyle \sum_{i_{k}=1}^{N_{k}} \tilde{F}(i_{1},i_{2}, \ldots{} , i_{k})\,,
\end{equation}
where $N_{j} = N + 2j$, $(j\in{\mathbb{N}})$, and $H(x)$ - is the unit step function (see Eq.(\ref{i33})).

Let us now define a global index $q$ as follows:
\begin{equation}\label{i37}
q = i_{1} + N_{1}(i_{2}-1) + N_{1}N_{2}(i_{3}-1) + \ldots{} + \Bigg[\displaystyle\prod_{j=1}^{k-1}N_{j}\Bigg](i_{k}-1)\,.
\end{equation}
Note, that the minimum value $q_{\text{min}}=1$ corresponds to the case then every index number is equal to unity. The maximum value $q_{\text{max}}=N_{1}\cdot{}N_{2}\cdot{}\,\cdots{}\,\cdot{}N_{k}$ corresponds to the case then every index number reaches its maximum, since the relation:
\begin{equation}\label{i38}
N_{1} + N_{1}(N_{2}-1) + N_{1}N_{2}(N_{3}-1) + \ldots{} + \Bigg[\displaystyle\prod_{j=1}^{k-1}N_{j}\Bigg](N_{k}-1) = \Bigg[\displaystyle\prod_{j=1}^{k}N_{j}\Bigg].
\end{equation} 
Now we construct one-to-one mapping between the global index $q$ and the index numbers $(i_{1},\,\ldots{}\,,i_{k})$. For that we need to solve the Eq.(\ref{i37}), i.e. to express any certain index number as a function of $q$ (not depending on any other index numbers). We do that separately for $i_{1}$, $i_{k}$, and all others $i_{r},\,(2\leq{}r\leq{}[k-1])$. \\
\underline{Obtaining of $i_{1}$}:\\
Since $(i_{1}-1)<N_{1}$, the expression $\big\{(q-1)/N_{1}\big\}$ does not depend on any index numbers except $i_{1}$. Thus, we obtain:
\begin{equation}\label{i39}
i_{1} = i_{1}(q) = 1 + N_{1}\big\{(q-1)/N_{1}\big\}.
\end{equation} 
\underline{Obtaining of $i_{k}$}:\\
Since $i_{1} + N_{1}(i_{2}-1) + N_{1}N_{2}(i_{3}-1) + \ldots{} + \Bigg[\displaystyle\prod_{j=1}^{k-2}N_{j}\Bigg](i_{k-1}-1) - 1 \, < \, \Bigg[\displaystyle\prod_{j=1}^{k-1}N_{j}\Bigg]$, (see Eq.(\ref{i37})), the expression $\Bigg\llbracket(q-1) \bigg / \,\displaystyle\prod_{j=1}^{k-1}N_{j}\Bigg\rrbracket$ does not depend on any index numbers except $i_{k}$. Thus we obtain:
\begin{equation}\label{i40}
i_{k} = i_{k}(q) = 1 + \Bigg\llbracket(q-1) \bigg / \,\displaystyle\prod_{j=1}^{k-1}N_{j}\Bigg\rrbracket \,.
\end{equation}
\underline{Obtaining of $i_{r},\,(2\leq{}r\leq{}[k-1])$}:\\
Here we make two steps, similar to the previous cases, in order to eliminate the dependence of higher and lower index numbers separately. Consider the following expression:
\begin{equation}\label{i41}
h_{r} = h_{r}(i_{1},\,\ldots\,,i_{r}) = i_{1} + \ldots{} + \Bigg[\displaystyle\prod_{j=1}^{r-1}N_{j}\Bigg](i_{r}-1)\,.
\end{equation}
According to Eq.(\ref{i38}) we have: $h_{r}-1\,<\,\Bigg[\displaystyle\prod_{j=1}^{r}N_{j}\Bigg]$. In addition, we have simple relations: $\big\{u/v\big\}\cdot{}v\equiv{u}\,$ for $\,u<v\,$, and $\big\{u/v\big\}\cdot{}v\equiv{0}\,$ for $\,(u/v)\in{\mathbb{N}}\,$. Thus we obtain:
\begin{equation}\label{i42}
h_{r} = 1 + \Bigg\{(q-1) \bigg / \,\displaystyle\prod_{j=1}^{r}N_{j}\Bigg\}\cdot{}\displaystyle\prod_{j=1}^{r}N_{j}\,\,.
\end{equation}
Next, we perform a step similar to the one we performed when obtain $i_{k}$ (see also Eq.(\ref{i40})).
\begin{equation}\label{i43}
i_{r} = 1 + \Bigg\llbracket(h_{r}-1) \bigg / \,\displaystyle\prod_{j=1}^{r-1}N_{j}\Bigg\rrbracket \,.
\end{equation}
Finally, we could write:
\begin{equation}\label{i44}
i_{r} = i_{r}(q) = 1 + \Bigg\llbracket{}N_{r}\Bigg\{(q-1) \bigg / \,\displaystyle\prod_{j=1}^{r}N_{j}\Bigg\}\Bigg\rrbracket \,.
\end{equation}
As can be seen, Eq.(\ref{i39}) corresponds to Eq.(\ref{i44}) if we put formally $r=1$. In addition, we could write \newline $(q-1)\,\equiv{} \Bigg\{(q-1) \bigg / \,\displaystyle\prod_{j=1}^{k}N_{j}\Bigg\}\cdot{}\displaystyle\prod_{j=1}^{k}N_{j}\,$. This means that Eq.(\ref{i40}) also corresponds to Eq.(\ref{i44}) if we put formally $r=k$. Thus we can write the described substitution for the \textit{discrete dimensional-convolution procedure} in the following form:
\begin{equation}\label{i45}
\begin{aligned}
&q = i_{1} + N_{1}(i_{2}-1) + N_{1}N_{2}(i_{3}-1) + \ldots{} + \Bigg[\displaystyle\prod_{j=1}^{k-1}N_{j}\Bigg](i_{k}-1)\,, &&\,q\in{}\Bigg[1,\,\ldots\,,\displaystyle\prod_{j=1}^{k}N_{j}\Bigg]\,;\\
&i_{r} = 1 + \Bigg\llbracket{}N_{r}\Bigg\{(q-1) \bigg / \,\displaystyle\prod_{j=1}^{r}N_{j}\Bigg\}\Bigg\rrbracket\,, &&\,r\in{}\Bigg[1,\,\ldots\,,k\Bigg]\,.
\end{aligned}
\end{equation} 
According to the previous equations every certain unique index number combination corresponds to the unique value of $q$. In addition, the total number of different index number combinations is equal to the total number of different values of $q$. Thus, we can conclude that Eq.(\ref{i45}) describes the one-to-one mapping for Eq.(\ref{i35}), and we can write:
\begin{equation}\label{i46}
E_{k} = \displaystyle \sum_{q=1}^{N_{1}\cdot{}N_{2}\cdot{}\,\cdots\,\cdot{}N_{k}} \; F\Big(i_{1}(q), \ldots{} , i_{k}(q)\Big) \,.
\end{equation}
Summarizing, the proposed procedure allows one to calculate a reduced R-sum without recursive indexing.

\section{Closed-form solution of the canonical three-term recurrence relation}

In this section we apply the previous results, in order to obtain the closed-form solution of Eq.(\ref{i4}) without recursive indexing. By using the exact solution via R-sum expansion Eq.(\ref{i21}), and discrete dimensional-convolution procedure Eq.(\ref{i36}), Eq.(\ref{i45}-\ref{i46}), we finally obtain:

\begin{equation}\label{i47}
\begin{aligned}
\vspace{0.2cm}
&a_{n+1} = 1 + \displaystyle \sum_{p=1}^{\llbracket{}(n+1)/2\rrbracket{}} \; \displaystyle \sum_{q=1}^{M(n,p)} G(n,p,q) \,\,;  &&&&&&&&&&\\
\vspace{0.3cm}
&M(n,p) = \displaystyle\prod_{r=1}^{p}(n-2p+2r) = \displaystyle\prod_{l=0}^{p-1}(n-2l)\,,\\
\vspace{0.2cm}
&G(n,p,q) = d_{g(1,n,p,q)} \cdot{} \displaystyle\prod_{m=2}^{p} \Bigg[d_{g(m,n,p,q)}\cdot{}H\Big[g(m,n,p,q)-g(m-1,n,p,q)-2\Big]\Bigg] \,,\\
\vspace{0.2cm}
&g(1,n,p,q) = q - (n-2p+2) \Bigg\llbracket(q-1) \bigg / \,(n-2p+2)\Bigg\rrbracket \,,\\
\vspace{0.2cm}
&g(m,n,p,q) = 1 + \Bigg\llbracket(q-1) \bigg / \,\displaystyle\prod_{j=1}^{m-1}(n-2p+2j)\Bigg\rrbracket - (n-2p+2m) \Bigg\llbracket(q-1) \bigg / \,\displaystyle\prod_{j=1}^{m}(n-2p+2j)\Bigg\rrbracket \,.
\end{aligned}
\end{equation} 
\begin{flushleft}\underline{Remark}:\end{flushleft} For the definition of the integer-valued function $g$ we use Eq.(\ref{i45}) and Gauss notation for the floor function: \newline $\llbracket{x}\rrbracket=\left\lfloor x \right\rfloor= x-\{x\} =\max\big(m\in{\mathbb{Z}} \; | \; m\leq{}x\big)$; \, \\ 
\vspace{0.0cm}\\
$H(x)$ is the unit step function, see Eq.(\ref{i33}).

\newpage

\section{Conclusions}

	In summary, we obtain the the closed-form solution of a canonical three-term recurrence relation Eq.(\ref{i4}) with an arbitrary given $n$-dependent coefficient $d_{n}$. The final non-recursive expression Eq.(\ref{i47}) includes a finite number of elementary operations and functions.
	
	Possible applications of the developed approaches, namely the R-sum theory and the discrete dimensional-convolution procedure, are not limited by the considered statement. Due to its universality, they could be used for solving other recursive problems, in particular many-term recurrence relations.
	
	An interesting and open question for the author, is how the solution Eq.(\ref{i47}) could be efficiently used for approximations and solving differential equations.

\end{document}